\documentclass{tran-l}

\usepackage{amscd,amssymb}
\usepackage[mathcal]{euscript}

\newtheorem{theorem}{Theorem}[section]
\newtheorem{lemma}[theorem]{Lemma}
\newtheorem{proposition}[theorem]{Proposition}

\newtheorem{corollary}[theorem]{Corollary}
\theoremstyle{definition}
\newtheorem{definition}[theorem]{Definition}

\newtheorem{example}[theorem]{Example}

\numberwithin{equation}{subsection}

\DeclareMathOperator{\ho}{Ho}
\DeclareMathOperator{\Hom}{Hom}

\newcommand{\boxprod}{\mathbin\square}
\newcommand{\cat}[1]{\mathcal{#1}}
\newcommand{\cof}{\text{-cof}}
 
\newcommand{\cell}{\text{-cell}} 
\newcommand{\inj}{\text{-inj}}
\renewcommand{\mod}{\text{-mod}}
\newcommand{\alg}{\text{-alg}}
\newcommand{\llp}{left lifting property with respect to }

\newcommand{\rlp}{right lifting property with respect to }
\newcommand{\sm}{\wedge}

\newcommand{\colim}{\operatorname{colim}}

\newcommand{\monoids}[1]{\mathcal{M}(\mathcal{#1})}

\newcommand{\ie}{\textit{i.e., }}
\newcommand{\mathcolon}{\colon\,}
\newcommand{\uc}{\textup{:}}
\newcommand{\ulp}{\textup{(}}
\newcommand{\urp}{\textup{)}}
\begin{document}

\title{Monoidal model categories}

\date{\today} 
\author{Mark Hovey} 
\address{Department of Mathematics \\ 
Wesleyan University \\ 
Middletown, CT}
\email{hovey@member.ams.org} 


\begin{abstract}
A monoidal model category is a model category with a closed monoidal
structure which is compatible with the model structure.  Given a
monoidal model category, we consider the homotopy theory of modules over
a given monoid and the homotopy theory of monoids.  We make minimal
assumptions on our model categories; our results therefore are more
general, yet weaker, than the results
of~\cite{schwede-shipley-monoids}.  In particular, our results apply to
the monoidal model category of topological symmetric
spectra~\cite{hovey-shipley-smith}.  
\end{abstract}

\maketitle

\section*{Introduction}

A monoidal model category is a (closed) monoidal category that is also a
model category in a compatible way.  Monoidal model categories abound in
nature: examples include simplicial sets, compactly generated
topological spaces, and chain complexes of modules over a commutative
ring.  The thirty-year long search for a monoidal model category of
spectra met success with the category of $S$-modules
of~\cite{elmendorf-kriz-mandell-may} and the symmetric spectra
of~\cite{hovey-shipley-smith}.

Given any monoidal category, one has categories of monoids and of
modules over a given monoid.  If we are working in a monoidal model
category, we would like these associated categories also to be model
categories, so that we can have a homotopy theory of rings and modules.
The first results on this subject were obtained
in~\cite{schwede-shipley-monoids}.  This paper is a followup to that
paper.  In~\cite{schwede-shipley-monoids}, the authors added the
following three assumptions about a monoidal model category $\cat{C}$:

\begin{enumerate}
\item [(a)] Every object of $\cat{C}$ is small relative to the whole
category; 
\item [(b)] $\cat{C}$ satisfies the monoid axiom; and 
\item [(c)] Given a monoid $A$ and a cofibrant left $A$-module $M$,
smashing over $A$ with $M$ takes weak equivalences of right $A$-modules
to weak equivalences.
\end{enumerate}

The first two assumptions guarantee the existence of a model structure
on the category of monoids and on the category of modules over a given
monoid.  The third assumption guarantees that a weak equivalence of
monoids induces a Quillen equivalence of the corresponding module
categories.  

All these assumptions are reasonable ones in any combinatorial
situation, such as simplicial sets, chain complexes, or simplicial
symmetric spectra.  However, for any category of topological spaces the
third assumption will fail, and the first assumption is not known to be
true and probably fails.  Furthermore, in the category of topological
symmetric spectra the second assumption is not known to hold.  

The goal of this paper, then, is to investigate what can be said when
these assumptions do not hold.  After a preliminary section reminding
the reader of some basic definitions and facts about model categories,
we begin in the second section by showing that one always gets a model
category of modules over a cofibrant monoid.  Furthermore, under a minor
assumption on our model category $\cat{C}$, we show that a weak
equivalence of cofibrant monoids induces a Quillen equivalence of the
corresponding module categories.  Also, a Quillen equivalence $F$ of
monoidal model categories induces a Quillen equivalence between
$R$-modules and $FR$-modules, for $R$ a monoid which is cofibrant in the
domain of $F$.  In the third section, we show that, if
the unit $S$ in $\cat{C}$ is cofibrant, then, though we do not get a
model category of monoids in general, we do at least get a homotopy
category of monoids.  In particular, given a general monoid $A$, we can
find a cofibrant monoid $QA$ and a weak equivalence and homomorphism
$QA\xrightarrow{}A$.  Then the model category $QA\mod $ is the homotopy
invariant replacement for the category $A\mod $, which may not even be a
model category.  We also show that the homotopy category of monoids is
itself homotopy invariant.  In particular, there is a homotopy category
of monoids of topological symmetric spectra, and this homotopy category
is equivalent to the homotopy category of monoids of simplicial
symmetric spectra.  

The most obvious question left unaddressed in this paper concerns the
category of \emph{commutative} monoids in a symmetric monoidal model
category.  What do we need to know to get a model structure on
commutative monoids?  Can we get a homotopy category of commutative
monoids in any symmetric monoidal model category?  The author does not
know the answer to these questions. 

The author would like to thank his coauthors Brooke Shipley and Jeff
Smith.  This paper grew out of~\cite{hovey-shipley-smith}, when the
authors of that paper realized that topological spaces are not as simple
as they had originally thought.  The author would also like to thank
Gaunce Lewis and Peter May for helping him come to that realization,
which of course they have understood for years.  

\section{Basics}

We will have to assume some familiarity with model categories on the
part of the reader.  A gentle introduction to the subject can be found
in~\cite{dwyer-spalinski}.  A more thorough and highly recommended
source is~\cite[Part 2]{hirschhorn}.  Other sources
include~\cite{hovey-model} and~\cite{kan-model}.

In particular, in a model category $\cat{C}$, we have a cofibrant
replacement functor $Q$ and a fibrant replacement functor $R$.  There
is a natural trivial fibration $QX\xrightarrow{q}X$, and $QX$ is
cofibrant.  Similarly, there is a natural trivial cofibration
$X\xrightarrow{}RX$ and $RX$ is fibrant.

Our basic object of study is a monoidal model category, which we now
define.  In a monoidal category $\cat{C}$, we will denote the monoidal
product by $\sm $ and the unit by $S$.  Note that in model category
theory, functions seem to come in adjoint pairs.  We will therefore
consider a \emph{closed} monoidal category rather than a general 
monoidal category.  This means that both functors $X\sm -$ and
$-\sm X$ have right adjoints natural in $X$.  For our purposes, the
closed structure just guarantees for us that both functors $X\sm -$ and
$-\sm X$ preserve colimits.  

\begin{definition}\label{defn-pushout-smash} 
Suppose $\cat{C}$ is a closed monoidal category.  Given maps
$f\mathcolon A\xrightarrow{}B$ and $g\mathcolon X\xrightarrow{}Y$ in
$\cat{C}$, define the \emph{pushout smash product} $f\boxprod g$ of $f$
and $g$ to be the map $(A\sm Y)\amalg _{A\sm X}(B\sm X)\xrightarrow{}B\sm
Y$.
\end{definition}

\begin{definition}\label{defn-monoidal-model}
Suppose $\cat{C}$ is a closed monoidal category which is also a model
category.  Then $\cat{C}$ is a \emph{monoidal model category} if the
following conditions hold.
\begin{enumerate}
\item [(a)] If $f$ and $g$ are cofibrations, so is $f\boxprod g$.  If
one of $f$ or $g$ is in addition a weak equivalence, so is $f\boxprod
g$.  
\item [(b)] Both maps $q\sm X\mathcolon QS\sm X\xrightarrow{}S\sm X\cong
X$ and $X\sm q\mathcolon X\sm QS\xrightarrow{}X$ are weak equivalences
for all cofibrant $X$.
\end{enumerate}
\end{definition}

The second condition is a consequence of the first in case the unit $S$
is cofibrant.  This is usually the case, but $S$ is not cofibrant in the
category of $S$-modules of~\cite{elmendorf-kriz-mandell-may}.  Without
the second condition, the homotopy category of a monoidal model category
would not be a monoidal category, because there would not be a unit.
With it, it is an exercise in derived functors, carried out
in~\cite[Chapter 4]{hovey-model}, to verify that the homotopy category
is indeed a monoidal category.

Of course, we only need half the second condition in case $\cat{C}$ is
symmetric monoidal, as it usually is in our examples.  

We point out, following the insight of Stefan Schwede, that the second
condition in Definition~\ref{defn-monoidal-model} is equivalent to
requiring that both maps $X\xrightarrow{}\Hom _{\ell }(QS,X)$ and
$X\xrightarrow{}\Hom _{r}(QS,X)$ are weak equivalences for all fibrant
$X$, where $\Hom _{\ell }$ and $\Hom _{r}$ are the two adjoints that
define the closed structure on $\cat{C}$.  To see this, one can show
that both the $\Hom $ conditions just defined and the $\sm $ conditions
of Definition~\ref{defn-monoidal-model} are equivalent to the unit
axioms in the monoidal category $\ho \cat{C}$.

Examples of symmetric monoidal model categories include the categories of
simplicial sets, compactly generated topological spaces,
$S$-modules~\cite{elmendorf-kriz-mandell-may}, symmetric
spectra~\cite{hovey-shipley-smith}, and topological symmetric spectra.  

The reader should note that, in a monodial model category, smashing with
a cofibrant object preserves cofibrations and trivial cofibrations, and
hence also weak equivalence between cofibrant objects, by Ken Brown's
lemma~\cite[Lemma 9.9]{dwyer-spalinski}.  

We will now repeat some standard definitions. 

\begin{definition}\label{defn-llp}
A map $f\mathcolon A\xrightarrow{}B$ in a category $\cat{C}$ is said to
have the \emph{\llp}another map $g\mathcolon X\xrightarrow{}Y$ if, for every
commutative square
\[
\begin{CD}
A @>>> X \\
@VfVV @VVgV \\
B @>>> Y
\end{CD}
\]
there is a lift $h\mathcolon B\xrightarrow{}X$ making the diagram
commute.  We also say that $g$ has the \emph{\rlp}$f$ in this situation. 
\end{definition}

The following argument is often used in model category theory. 

\begin{proposition}[The Retract Argument]\label{prop-retract}
Let $\cat{C}$ be a category and let $f=pi$ be a factorization in $\cat{C}$.
\begin{enumerate}
\item If $p$ has the \rlp $f$ then $f$ is a retract of $i$.
\item If $i$ has the \llp $f$ then $f$ is a retract of $p$.
\end{enumerate}
\end{proposition}

\begin{proof}
We only prove the first part, as the second is similar.  Since $p$
has the right lifting property with respect to $f$, we have a lift
$g\mathcolon Y\xrightarrow{}Z$ in the diagram
\[
\begin{CD}
X @>i>> Z \\
@VfVV @VpVV \\
Y @= Y
\end{CD}
\]
This gives a diagram
\[
\begin{CD}
X @= X @= X\\
@VfVV @ViVV @VfVV \\
Y @>g>> Z @>p>> Y
\end{CD}
\]
where the horizontal compositions are identity maps, showing that $f$ is
a retract of $i$.
\end{proof}

In model category theory, one very often needs to construct
factorizations.  The standard way to do this is by the small object
argument.  

\begin{definition}\label{defn-small}
Suppose $A$ is an object of a cocomplete category $\cat{C}$.  Suppose
$\cat{D}$ is a subcategory of $\cat{C}$.  We say that $A$ is
\emph{small} relative to $\cat{D}$ if there is a cofinal class $S$
of ordinals such that, for all $\alpha \in S$ and for all colimit-preserving
functors $X\mathcolon \alpha \xrightarrow{}\cat{C}$ such that each map
$X_{\beta }\xrightarrow{}X_{\beta +1}$ is in $\cat{D}$, the induced map 
\[
\colim_{\beta <\alpha } \cat{C}(A,X_{\beta }) \xrightarrow{}
\cat{C}(A,\colim _{\beta <\alpha }X_{\beta }) 
\]
is an isomorphism.  
\end{definition}

In this context, ``cofinal'' means that given any ordinal $\alpha $
there is an ordinal $\beta \in S$ with $\alpha \leq \beta $.  For
example, every set is small relative to the whole category of sets.  For
a finite set $A$, the class $S$ is the collection of limit ordinals; for
a more general set $A$ one has to take a sparser collection of ordinals.

\begin{definition}\label{defn-I-cell}
Suppose $I$ is a collection of maps in a cocomplete category $\cat{C}$.
Define \emph{$I\inj $} to be the class of all maps with the \rlp $I$, and
define \emph{$I\cof $} to be the class of all maps with the \llp $I\inj $.
Define \emph{$I\cell $} to be the class of all transfinite compositions of
pushouts of $I$.  That is, for any map $f\mathcolon A\xrightarrow{}B$ in
$I\cell $ there is an ordinal $\alpha $ and a colimit-preserving functor
$X\mathcolon \alpha \xrightarrow{}\cat{C}$ such that $X_{0}=A$, the map
$X_{0}\xrightarrow{}\colim_{\beta <\alpha } X_{\beta }$ is isomorphic
to $f$, and each map $X_{\beta }\xrightarrow{}X_{\beta +1}$ is a pushout
of a map of $I$.  
\end{definition}

Note that $I\cell \subseteq I\cof $.  

\begin{theorem}[The Small Object Argument]\label{thm-small-object}
Suppose $I$ is a set of maps in a cocomplete category $\cat{C}$, and
suppose that the domains of $I$ are small relative to $I\cell $.  Then
there is a functorial factorization of every map in $\cat{C}$ into a
map of $I\cell $ followed by a map of $I\inj $.  
\end{theorem}

We will not prove this theorem: see~\cite[Section 12.4]{hirschhorn}
or~\cite[Section 2.1]{hovey-model}.  Note that most authors include
transfinite compositions of pushouts of coproducts of $I$, but this is
not necessary in view of~\cite[Proposition 12.2.5]{hirschhorn}.  Also,
if the domains of $I$ are small relative to $I\cell $, then every map of
$I\cof $ is a retract of a map of $I\cell $ by the retract argument, and
furthermore the domains of $I$ are small relative to $I\cof
$~\cite[Theorem 12.4.21]{hirschhorn}.  

We can now define a cofibrantly generated model category.  

\begin{definition}\label{defn-cofib-gen}
A model category $\cat{C}$ is \emph{cofibrantly generated} if there is
are sets $I$ and $J$ of maps of $\cat{C}$ such that the following
conditions hold.  
\begin{enumerate}
\item The domains of $I$ are small relative to the cofibrations.  The
domains of $J$ are small relative to the trivial cofibrations. 
\item The fibrations form the class $J\inj $.  The trivial fibrations
form the class $I\inj $. 
\end{enumerate}
\end{definition}

Every model category in common use is cofibrantly generated, such as the
model categories of topological spaces, simplicial sets, and chain
complexes.  In a cofibrantly generated model category, the cofibrations
form the class $I\cof $ and the trivial cofibrations form the class
$J\cof $.  

The big advantage of cofibrantly generated model categories is that they
allow us to prove things by induction.  Suppose we have some claim about
cofibrant objects $A$ in a cofibrantly generated model category.  That
claim is almost certain to be preserved by retracts, so we can usually
assume the map $0\xrightarrow{}A$ is in $I\cell $, so is a transfinite
composition of pushouts of maps of $I$.  We can then use transfinite
induction.  

We also remind the reader that functors between model categories come in
adjoint pairs.  A functor $F\mathcolon \cat{C}\xrightarrow{}\cat{D}$
between model categories with right adjoint $U$ is called a \emph{left
Quillen functor} (and $U$ is called a \emph{right Quillen functor}) if
$F$ preserves cofibrations and trivial cofibrations.  Equivalently, we
can require that $U$ preserve fibrations and trivial fibrations.  A
Quillen pair induces a pair of adjoint functors $LF$ and $RU$ on the
homotopy categories.  The functor $LF$ is defined by $(LF)(X)=F(QX)$,
where $Q$ is the functorial cofibrant replacement functor (this is a
very good reason to assume these factorizations are part of the model
structure, as is done in~\cite{hovey-model}).  Simiilarly $RU$ is
defined by $(RU)(X)=U(RX)$, where $R$ is the functorial fibrant
replacement functor.  

A left Quillen functor $F$ is called a \emph{Quillen equivalence} if for
all cofibrant $A\in \cat{C}$ and fibrant $X\in \cat{D}$, a map
$FA\xrightarrow{}X$ is a weak equivalence if and only if its adjoint
$A\xrightarrow{}UX$ is a weak equivalence.

The following lemma deserves to be better known than it is.  Recall that
a functor is said to \emph{reflect} some property of morphisms if, given
a morphism $f$, if $Ff$ has the property so does $f$.  
 
\begin{lemma}\label{lem-Quillen-equiv}
Suppose $F\mathcolon \cat{C}\xrightarrow{}\cat{D}$ is a Quillen functor
with right adjoint $U$.  The following are equivalent\uc 
\begin{enumerate}
\item [(a)] $F$ is a Quillen equivalence.
\item [(b)] $F$ reflects weak equivalences between cofibrant
objects and, for every fibrant $Y$, the map $FQUY\xrightarrow{}Y$ is a
weak equivalence. 
\item [(c)] $U$ reflects weak equivalences between fibrant objects and,
for every cofibrant $X$, the map $X\xrightarrow{}URFX$ is a weak
equivalence.  
\item [(d)] $LF$ is an equivalence of categories. 
\end{enumerate}
\end{lemma}

\begin{proof}
Suppose first that $F$ is a Quillen equivalence.  Then, if $X$ is
cofibrant, the weak equivalence $FX\xrightarrow{}RFX$ gives rise to a
weak equivalence $X\xrightarrow{}URFX$.  Similarly, the weak equivalence
$QUX\xrightarrow{}UX$ gives rise to a weak equivalence
$FQUX\xrightarrow{}X$.  This shows that (a) implies half of (b) and
(c).  Now suppose $f\mathcolon X\xrightarrow{}Y$ is a map between
cofibrant objects such that $Ff$ is a weak equivalence.  Since both maps
$X\xrightarrow{}URFX$ and $Y\xrightarrow{}URFY$ are weak equivalences,
$f$ is a weak equivalence if and only if $URFf$ is a weak equivalence.
Since $Ff$ is a weak equivalence, $R$ preserves weak equivalences, and
$U$ preserves weak equivalences between fibrant objects, we find that
$f$ is a weak equivalence.  Thus (a) implies (b), and a similar argument
shows that (a) implies (c).  

To see that (b) implies (d), note that the counit map
$(LF)(RU)X\xrightarrow{}X$ is an isomorphism by hypothesis.  We must
show that the unit map $X\xrightarrow{}(RU)(LF)X$ is an isomorphism.
But $(LF)X\xrightarrow{}(LF)(RU)(LF)X$ is inverse to the counit map of
$(LF)X$, so is an isomorphism.  Since $F$ reflects weak equivalences
between cofibrant objects, this implies that $QX\xrightarrow{}QURFQX$ is
a weak equivalence.  Since $Q$ reflects all weak equivalencs, this
implies that $X\xrightarrow{}URFQX=(RU)(LF)X$ is a weak equivalence, as
required.  A similar proof shows that (c) implies (d).  

To see that (d) implies (a), note that $(LF)X$ is isomorphic to $FX$ in
the homotopy category when $X$ is cofibrant, and similarly $(RU)Y$ is
isomorphic to $UY$ in the homotopy category when $U$ is fibrant.  So
$FX\xrightarrow{}Y$ is a weak equivalence if and only
$(LF)X\xrightarrow{}Y$ is an isomorphism in the homotopy category.
Since $LF$ is an equivalence of categories with adjoint $RU$, this is
true if and only if $X\xrightarrow{}(RU)Y$ is an isomorphism.  But this
holds if and only if $X\xrightarrow{}UY$ is a weak equivalence, as
required.  
\end{proof}

\section{Modules}

In this section we investigate model categories of modules over a monoid
$A$ in a cofibrantly generated monoidal model category $\cat{C}$.  

\begin{theorem}\label{thm-cofibrant-modules}
Suppose $\cat{C}$ is a cofibrantly generated monoidal model category
with generating cofibrations $I$ and generating trivial cofibrations
$J$.  Let $A$ be a monoid in $\cat{C}$, and suppose the following
conditions hold.
\begin{enumerate}
\item The domains of $I$ are small relative to $(A\sm I)\cell $. 
\item The domains of $J$ are small relative to $(A\sm J)\cell $. 
\item Every map of $(A\sm J)\cell $ is a weak equivalence. 
\end{enumerate}
Then there is a cofibrantly generated model structure on the category of
left $A$-modules, where a map is a weak equivalence or fibration if and
only if it is a weak equivalence or fibration in $\cat{C}$.
\end{theorem}

There is an obvious analogue of this theorem for right $A$-modules. 

\begin{proof}
By adjointness, the fibrations of $A$-modules form the class $(A\sm
J)\inj $, and the trivial fibrations form the class $(A\sm I)\inj $.  We
therefore define the cofibrations of $A$-modules to be the class $(A\sm
I)\cof $, and take our generating trivial cofibrations to be $A\sm J$.
Since each element of $J$ is in $I\cof $, each element of $A\sm J$ is in
$(A\sm I)\cof $, and so the maps of $(A\sm J)\cof $ are cofibrations in
$A\mod $.

The category of $A$-modules is certainly bicomplete, with limits and
colimits taken in $\cat{C}$.  The retract and two out of three axioms
are immediate, as is the lifting axiom for cofibrations and trivial
fibrations.  By assumption, the domains of $I$ are small relative to
$(A\sm I)\cell $.  By adjointness, it follows that the domains of $A\sm
I$ are small in $A\mod $ relative to $(A\sm I)\cell $.  Thus the small
object argument gives us the cofibration-trivial fibration half of the
factorization axiom.  Similarly, the domains of $A\sm J$ are small in
$A\mod $ relative to $(A\sm J)\cell $.  We can then factor any map in
$A\mod $ into a map of $(A\sm J)\cell $ followed by a fibration.  We
have already seen that the maps of $(A\sm J)\cell $ are cofibrations,
and by assumption they are weak equivalences.  This gives the
other half of the factorization axiom.

For the remaining lifting axiom, suppose $f$ is a cofibration
and weak equivalence.  Factor $f=pi$, where $i\in (A\sm J)\cell $ and $p$
is a fibration of $A$-modules.  Since $i$ is a weak equivalence, and so
is $f$, it follows that $p$ is a weak equivalence.  Thus $f$ has the
\llp $p$.  By the Retract Argument~\ref{prop-retract}, $f$ is a retract
of $i$, and so has the \llp all fibrations of $A$-modules, as required.
\end{proof}

\begin{corollary}\label{cor-cofibrant-modules}
Suppose $\cat{C}$ is a cofibrantly generated monoidal model category,
and suppose $A$ is a monoid which is cofibrant in $\cat{C}$.  Then there
is a cofibrantly generated model structure on the category of
\ulp left, or right\urp{} $A$-modules, where a map is a weak
equivalence or fibration if and only if it is a weak equivalence or
fibration in $\cat{C}$.  Furthermore, a cofibration of $A$-modules is a
cofibration in $\cat{C}$.  
\end{corollary}

\begin{proof}
For definiteness, we work with left $A$-modules.  Since $A$ is cofibrant
in $\cat{C}$, every map of $A\sm I$ is a cofibration, so $(A\sm I)\cell
\subseteq I\cof $.  Since the domains of $I$ are small relative to
$I\cof $, they are certainly small relative to $(A\sm I)\cell $.
Similarly, $(A\sm J)\cell \subseteq J\cof $, so the domains of $J$ are
small relative to $(A\sm J)\cell $, and the maps of $(A\sm J)\cell $ are
weak equivalences.
\end{proof}

As an example, we would like to consider the model category of
topological spaces.  This is not, however, a monoidal model category,
since the functor $X\times -$ need not have a right adjoint unless $X$
is locally compact Hausdorff.  To get around this, it is usual to
consider some version of compactly generated spaces.  We use the
definitions of~\cite[Appendix A]{lewis-thesis}, the properties of which
are summarized in~\cite[Section 6.1]{hovey-shipley-smith}.  In
particular, a subset $U$ of a topological space $X$ is \emph{compactly open}
if, for every continuous $f\mathcolon K\xrightarrow{}X$, where $K$ is
compact Hausdorff, the preimage $f^{-1}(U)$ is open.  The space $X$ is
called a \emph{$k$-space} if every compactly open space is open.  A
space is called \emph{compactly generated} if it is both a $k$-space and
weak Hausdorff; \ie for every continuous $f\mathcolon K\xrightarrow{}X$
where $K$ is compact Hausdorff, the image $f(K)$ is closed.  Then both
the category $\cat{K}$ of $k$-spaces and the category $\cat{T}$ of
compactly generated spaces are cofibrantly generated symmetric monoidal
model categories.  

In order to understand these categories a little better, we need the
following lemma. 

\begin{lemma}\label{lem-monoid-axiom-fibrant}
Suppose $\cat{C}$ is a cofibrantly generated monoidal model category,
such that the domains and codomains of the generating trivial
cofibrations $J$ are fibrant.  Suppose in addition that there is an
object $I\in \cat{C}$ and a factorization $S\amalg
S\xrightarrow{(i_{0},i_{1})}I\xrightarrow{}S$ of the fold map, where
$(i_{0},i_{1})$ is a cofibration, such that the induced map $X\sm
I\xrightarrow{}X$ is a weak equivalence for all $X$.  Then $\cat{C}$
satisfies the monoid axiom\uc{} that is, every map of $(\cat{C}\sm J)\cell
$ is a weak equivalence.
\end{lemma}

\begin{proof}
We only sketch the proof.  We can define a map $f\mathcolon
A\xrightarrow{}B$ to be a \emph{strong deformation retract} if there is
a retraction $r\mathcolon B\xrightarrow{}A$ such that $rf=1_{A}$ and a
homotopy $H\mathcolon B\sm I\xrightarrow{}B$ such that $Hi_{0}=rf$ and
$Hi_{1}=1_{B}$.  Here we are using the specific object $I$ in the
hypothesis of the lemma.  In particular, $B\sm I$ will not be a
cylinder object for $B$ in general.  Nevertheless, one can check that
any strong deformation retract us a weak equivalence, and furthermore,
that the class of strong deformation retracts is closed under smashing
with an arbitrary object, pushouts, and transfinite compositions.
Furthermore, following the argument of~\cite[pg. 2.5]{quillen-htpy}, we
can see that each map of $J$ is a strong deformation retract.  Thus
every map of $(\cat{C}\sm J)\cell $ is a strong deformation retract, and
hence a weak equivalence, as required.  
\end{proof}

In particular, this applies to both the category $\cat{K}$ of $k$-spaces
and the category $\cat{T}$ of compactly generated topological spaces,
where in this case $I$ is the usual unit interval.  So both $\cat{K}$
and $\cat{T}$ satisfy the monoid axiom.  However, the smallness
conditions of Theorem~\ref{thm-cofibrant-modules} do not appear to hold
in $\cat{K}$.  In any topological category, the best one can usually do
for smallness is that every object is small relative to the inclusions.
But inclusions are not preserved very well, so we do not know that the
maps of $(A\sm I)\cell $ are inclusions.  In $\cat{T}$, however,
\emph{closed} inclusions are preserved by almost every construction one
ever makes (see~\cite[Section 6.1]{hovey-shipley-smith}, based
on~\cite[Appendix A]{lewis-thesis}).  In particular, the maps of $(A\sm
I)\cell $ and of $(A\sm J)\cell $ are closed inclusions.  Therefore, we
do get model categories of modules over an arbitrary monoid in
$\cat{T}$.

In order for the model category of $A$-modules to be useful, it must be
both homotopy invariant in appropriate senses and have good properties.
We begin with the homotopy invariance.

\begin{theorem}\label{thm-maps-cofibrant-monoids}
Suppose $\cat{C}$ is a cofibrantly generated monoidal model category
such that the domains of the generating cofibrations can be taken to be
cofibrant.  Suppose $f\mathcolon A\xrightarrow{}A'$ is a weak
equivalence of monoids which are cofibrant in $\cat{C}$.  Then the
induction functor induced by $f$ and its right
adjoint, the restriction functor, define a Quillen equivalence from the
model category of left \ulp resp. right\urp{} $A$-modules to
the model category of left \ulp resp. right\urp{} $A'$-modules.
\end{theorem}

\begin{proof}
Again, we work with left modules for definiteness.  In this case the
induction functor takes $M$ to $A'\sm _{A}M$.  The restriction
functor obviously preserves weak equivalences and fibrations, so is a
right Quillen functor.  Furthermore, the restriction functor reflects
weak equivalences as well.  It follows from
Lemma~\ref{lem-Quillen-equiv} that induction is a Quillen equivalence if
and only if for all cofibrant $A$-modules $M$, the map
$M\xrightarrow{}A'\sm _{A}M$ is a weak equivalence.  Because the
category of $A$-modules is cofibrantly generated, we may as well assume
that $M$ is the colimit of a colimit-preserving functor $\alpha
\xrightarrow{}A\mod $, where $\alpha $ is an ordinal, $M_{0}=0$
and each map $M_{\beta }\xrightarrow{}M_{\beta +1}$ is a pushout of a
map of $A\sm I$.

We will prove by transfinite induction that the map $i_{\beta
}\mathcolon M_{\beta }\xrightarrow{}A'\sm _{A}M_{\beta }$ is a weak
equivalence for all $\beta \leq \alpha $, taking $M_{\alpha }=\colim
_{\alpha <\beta }M_{\beta }=M$.  Getting the induction started is easy,
since $M_{0}=0$.  For the successor ordinal case, suppose that $i_{\beta
}$ is a weak equivalence.  We have a pushout diagram 
\[
\begin{CD}
A\sm K @>>> A\sm L \\
@VVV @VVV \\
M_{\beta } @>>> M_{\beta +1}
\end{CD}
\]
where $K\xrightarrow{}L$ is some map of $I$.  Both horizontal maps are
cofibrations in $\cat{C}$.  Furthemore, because $K$ and $L$ are by
assumption cofibrant in $\cat{C}$, each object in the diagram is
cofibrant in $\cat{C}$.  By applying the functor $A'\sm _{A}-$, we
get an analogous pushout diagram
\[
\begin{CD}
A'\sm K @>>> A'\sm L \\
@VVV @VVV \\
A'\sm _{A}M_{\beta } @>>> A'\sm _{A} M_{\beta +1}
\end{CD}
\]
where again the horizontal maps are cofibrations in $\cat{C}$, and each
object is cofibrant in $\cat{C}$.  There is a map from the first pushout
square to the second, which by the induction hypothesis is a weak
equivalence on the lower left square.  Since $K$ is cofibrant, smashing
with $K$ preserves trivial cofibrations in $\cat{C}$, and hence, by Ken
Brown's lemma~\cite[Lemma 9.9]{dwyer-spalinski}, preserves weak
equivalences between cofibrant objects of $\cat{C}$.  Thus the map of
pushout squares is also a weak equivalence on the upper left corner, and
by similar reasoning, on the upper right corner.  Dan Kan's cubes lemma
(see~\cite[Section 5.2]{hovey-model} or~\cite{kan-model}) then shows
that it is a weak equivalence on the lower right corner as well.  Thus
$i_{\beta +1}$ is a weak equivalence.

Now consider the limit ordinal case.  Here we assume $i_{\gamma }$ is a
weak equivalence for all $\gamma <\beta $.  We then have a map of
sequences 
\[
\begin{CD}
M_{0} @>>> M_{1} @>>>\dots @>>> M_{\gamma } @>>> \dots \\
@Vi_{0}VV @Vi_{1}VV @. @Vi_{\gamma }VV @. \\
A'\sm _{A}M_{0} @>>> A'\sm _{A}M_{1} @>>> \dots @>>> A'\sm _{A}M_{\gamma
} @>>> \dots 
\end{CD}
\]
where the vertical maps are all weak equivalences, and the horizontal
maps are cofibrations of cofibrant objects.  Then~\cite[Proposition
18.4.1]{hirschhorn} implies that $i_{\beta }$ is also a weak
equivalence, as required. 
\end{proof}

The author learned the following example from Neil Strickland. 

\begin{example}\label{ex-top-maps}
Let us consider the category $\cat{T}_{*}$ of compactly generated
pointed spaces.  This is a monoidal model category under the smash
product.  Let $A$ denote the nonnegative natural numbers together with
infinity, given the discrete topology.  With infinity as the base point,
$A$ is a monoid in $\cat{T}_{*}$, with unit $0$.  Let $A'$ be the same
set as $A$, but with the one-point compactification topology.  Then $A'$
is also a topological monoid, with basepoint at infinity and unit $0$,
and the identity $A\xrightarrow{f}A'$ is a homomorphism and weak
equivalence of monoids.  Note that $A'$ is not cofibrant as a
topological space, so Theorem~\ref{thm-maps-cofibrant-monoids} does not
apply.  In fact, the induction functor does not define a Quillen
equivalence from $A$-modules to $A'$-modules.  Indeed, take $M=A\sm
S^{1}$.  Then the map $M\xrightarrow{}A'\sm _{A}M$ is just the
suspension of $f$.  The suspension of $A$ is an infinite wedge of
circles, but the suspension of $A'$ is the Hawaiian earring.  In
particular, $f\sm S^{1}$ is not a weak equivalence, so induction can not
be a Quillen equivalence.
\end{example}

In light of this example, it seems to the author that one should avoid
considering modules over monoids which are not cofibrant in $\cat{C}$,
just we generally avoid suspending non-cofibrant topological spaces.  
There are some categories, however, where a weak equivalence of monoids
always induces a Quillen equivalence of the corresponding module
categories.  This is true in $S$-modules~\cite[Theorem
3.8]{elmendorf-kriz-mandell-may} and in simplicial symmetric
spectra~\cite[Theorem 5.5.9]{hovey-shipley-smith}.  See~\cite[Theorem
3.3]{schwede-shipley-monoids}. 

Theorem~\ref{thm-maps-cofibrant-monoids} shows that the model category
of $A$-modules is homotopy invariant under weak equivalences of
cofibrant objects in $\cat{C}$.  But we would also like the model
category of $A$-modules to be homotopy invariant under weak
equivalences of monoidal model categories $\cat{C}$.

For this to make sense, we need a notion of monoidal Quillen functor.  

\begin{definition}\label{defn-monoidal-Quillen}
Suppose $F\mathcolon \cat{C}\xrightarrow{}\cat{D}$ is a left Quillen
functor between monoidal model categories.  Then $F$ is a \emph{monoidal
Quillen functor} if $F$ is monoidal and the induced map
$F(QS)\xrightarrow{}FS\cong S$ is a weak equivalence.   
\end{definition}

This second condition is easy to overlook; it is essential in case the
unit $S$ of $\cat{C}$ is not cofibrant in order to be sure $LF$
is a monoidal functor.  

\begin{theorem}\label{thm-modules-category}
Suppose $F\mathcolon \cat{C}\xrightarrow{}\cat{D}$ is a monoidal Quillen
equivalence of cofibrantly generated monoidal model categories, with
right adjoint $U$.  Suppose that $A$ is a monoid in $\cat{C}$ which is
cofibrant.  Then $F$ induces a Quillen equivalence $F\mathcolon A\mod
\xrightarrow{}FA\mod $.  
\end{theorem}

\begin{proof}
Since $F$ is monoidal, $FA$ is a monoid in $\cat{D}$; since $F$
preserves cofibrant objects, $FA$ is cofibrant in $\cat{D}$.  Given an
$A$-module $M$, $FM$ is an $FA$-module with structure map $FA\sm FM\cong
F(A\sm M)\xrightarrow{}FM$.  Hence $F$ does define a functor
$F\mathcolon A\mod \xrightarrow{}FA\mod $.  Let $\eta \mathcolon
X\xrightarrow{}UFX$ denote the unit of the adjunction, and let
$\varepsilon \mathcolon FUX\xrightarrow{}X$ denote the counit.  Given an
$FA$-module $N$, $UN$ is an $A$-module; its structure map is given by
the composite
\[
A\sm UN \xrightarrow{\eta \sm 1}UFA\sm UN \xrightarrow{}U(FA\sm
N)\xrightarrow{}UN
\]
where the second map is adjoint to the composite
\[
F(UFA\sm UN) \xrightarrow{\cong } FUFA \sm FUN \xrightarrow{\varepsilon
\sm \varepsilon } FA\sm N .
\]
One can easily check that the resulting functor $U\mathcolon FA\mod
\xrightarrow{}A\mod $ is right adjoint to $F$.  The functor $U$ clearly
preserves fibrations and trivial fibrations, so is a right Quillen
functor.  Furthermore, $U$ reflects weak equivalences between fibrant
objects, since $U$ does so as a functor from $\cat{D}$ to $\cat{C}$.  We
are of course using Lemma~\ref{lem-Quillen-equiv}, which tells us that
we need only check that the map $X\xrightarrow{}ULFX$ is a weak
equivalence for cofibrant $A$-modules $X$.  Here $L$ is a fibrant
replacement functor in $FA\mod $, not the fibrant replacement functor
$R$ in $\cat{D}$.  Nevertheless, there is a weak equivalence 
$FX\xrightarrow{}LFX$, and since $LFX$ is fibrant in $\cat{D}$, there is
a map $RFX\xrightarrow{}LFX$ in $\cat{D}$ which is necessarily a weak
equivalence.  Thus the map $URFX\xrightarrow{}ULFX$ is a weak
equivalence in $\cat{C}$.  Since $X$ is a cofibrant $A$-module, and in
particular cofibrant in $\cat{C}$, the map $X\xrightarrow{}URFX$ is a
weak equivalence.  It follows that the map $X\xrightarrow{}ULFX$ is a
weak equivalence as desired.  
\end{proof}

We now discuss some of the properties of the model category $A\mod $.  
If $\cat{C}$ is a symmetric monoidal category, and $A$ is a commutative
monoid, then it is well known that $A\mod $ is also a symmetric monoidal
category.  We would like the model structure on $A\mod $ to be
compatible with this symmetric monoidal structure.  

\begin{proposition}\label{prop-boxprod-R}
Suppose $\cat{C}$ is a cofibrantly generated symmetric monoidal model
category.  In addition, suppose that either
\begin{enumerate}
\item The unit $S$ is cofibrant and $A$ is a commutative monoid satisfying the
conditions of Theorem~\ref{thm-cofibrant-modules}; or 
\item $A$ is a commutative monoid cofibrant in $\cat{C}$.
\end{enumerate}
Then the model category $A\mod $ is a cofibrantly generated symmetric
monoidal model category.
\end{proposition}

\begin{proof}
It is well-known that $A\mod $ is closed symmetric monoidal:
see~\cite[Section 2.2]{hovey-shipley-smith} for details.  The symmetric
monoidal structure is denoted $\sm _{A}$, and we then have an analogous
definition of $\boxprod _{A}$.  Let $I$ be the set of generating
cofibrations of $\cat{C}$ and let $J$ be the set of generating trivial
cofibrations.  Then $A\sm I$ is the set of generating cofibrations of
$A\mod $, and $A\sm J$ is the set of generating trivial cofibrations.
We have
\[
(A\sm I)\boxprod _{A}(A\sm I)=(A\sm I)\boxprod I=A\sm (I\boxprod
I) \subseteq A\sm (I\cof )\subseteq (A\sm I)\cof .
\]
Thus $f\boxprod _{A}g$ is a cofibration of $A$-modules if $f$ and $g$
are, using~\cite[Lemma 2.3]{schwede-shipley-monoids}.  See
also~\cite[Corollary 5.3.5]{hovey-shipley-smith}.  A similar argument
shows that $f\boxprod _{A}g$ is a trivial cofibration if either $f$ or
$g$ is.  Thus the pushout product part of the definition of a monoidal
model category holds.  In the first case, we are done; since $S$ is
cofibrant in $\cat{C}$, $A$ is cofibrant in $A\mod $.  In the second
case, let $QS$ be a cofibrant replacement for the unit $S$ in $\cat{C}$,
so that $QS$ is cofibrant and we have a weak equivalence
$QS\xrightarrow{}S$.  Then $A\sm QS$ is cofibrant in $A\mod $, and,
since $A$ is cofibrant and $\cat{C}$ is monoidal, the $A$-module map
$A\sm QS\xrightarrow{}A$ is a weak equivalence.  Thus $A\sm QS$ is a
cofibrant replacement for the unit $A$ in $A\mod $.  We must show that,
if $M$ is a cofibrant $A$-module, then the map $(A\sm QS)\sm _{A}
M\xrightarrow{}A\sm _{A}M=M$ is still a weak equivalence.  But $(A\sm
QS)\sm _{A}M\cong QS\sm M$.  Since $M$ is cofibrant in $\cat{C}$, the
desired result holds.
\end{proof}

Note that the Quillen equivalences of
Theorems~\ref{thm-maps-cofibrant-monoids} and~\ref{thm-modules-category}
are monoidal Quillen equivalences in case the monoids involved are
commutative.

In case $A$ is not commutative, there is still a closed action of
$\cat{C}$ on $A\mod $; $M\sm X$ is an $A$-module if $M$ is an $A$-module
and $X$ is arbitrary.  This action also respects the model structures:
$f\boxprod g$ is a cofibration of $A$-modules if $f$ is a cofibration of
$A$-modules and $g$ is a cofibration.  Furthermore, $f\boxprod g$ is a
weak equivalence if either $f$ or $g$ is.  However, we will have trouble
with the unit unless we assume either $A$ or $S$ is cofibrant, just as
above.  

Since Hirschhorn's landmark treatment~\cite{hirschhorn}, it has become
clear that the right collection of model categories to work with is the
collection of left proper cellular model categories.  Hirschhorn shows
that one can perform Bousfield localization in this setting.  

\begin{proposition}\label{prop-cellular}
Suppose $\cat{C}$ is a left proper cellular monoidal model category.
Suppose $A$ is a monoid which is cofibrant in $\cat{C}$.  Then the model
category $A\mod $ is also left proper and cellular.
\end{proposition}

Because the definition of cellular is technical, the proof of this
proposition would take us too far afield.  It can be proved by the
methods of~\cite[Section 6]{hovey-stable-model}.  

\section{Algebras}

In this section we study the category $A\alg $ of algebras over a
commutative monoid $A$ in a cofibrantly generated symmetric monoidal
model category $\cat{C}$.  Note that in case $A=S$, an $S$-algebra is
the same thing as a monoid in $\cat{C}$.  Furthermore, an $A$-algebra is
just a monoid in the symmetric monoidal category $A\mod $.  We use this
to reduce to the case of monoids.

The obvious definitions to make for a model structure on $A\alg $ are
the following.  We define a homomorphism $f\mathcolon X\xrightarrow{}Y$
of $A$-algebras to be a \emph{weak equivalence} (\emph{fibration}) if
and only if $f$ is a weak equivalence (fibration) in $\cat{C}$.  Then
define $f$ to be a \emph{cofibration} if and only if $f$ has the \llp
all homomorphisms of $A$-algebras which are both weak equivalences and
fibrations.

Note that the forgetful functor $A\alg \xrightarrow{}\cat{C}$ has a left
adjoint, the free algebra functor $T$.  Of course, $T(X)=A\sm \coprod _{n\geq
0}X^{\sm n}$, where $X^{\sm 0}=S$.  We will usually use this only when
$A=S$.  

We begin by slightly generalizing the main result
of~\cite{schwede-shipley-monoids}.  

\begin{theorem}\label{thm-monoids-yes}
Suppose $\cat{C}$ is a cofibrantly generated symmetric monoidal model
category, with generating cofibrations $I$ and generating trivial
cofibrations $J$.  Suppose the following conditions hold.
\begin{enumerate}
\item The domains of $I$ are small relative to $(\cat{C}\sm I)\cell $. 
\item The domains of $J$ are small relative to $(\cat{C}\sm J)\cell $. 
\item The maps of $(\cat{C}\sm J)\cell $ are weak equivalences; \ie the
monoid axiom holds. 
\end{enumerate}
Let $A$ be a commutative monoid in $\cat{C}$.  Then $A\alg $ is a
cofibrantly generated model category where a map of $A$-algebras is a
weak equivalence or fibration if and only if it is so in $\cat{C}$.
\end{theorem}

\begin{proof}
We first point out that we can assume that $A=S$.  Indeed, our
conditions guarantee that $A\mod $ is a cofibrantly generated model
category (see Theorem~\ref{thm-cofibrant-modules}) which is symmetric
monoidal, and that the pushout product half of the definition of a
monoidal model category holds (see the proof of
Proposition~\ref{prop-boxprod-R}).  Furthermore, we have $A\mod \sm
_{A}(A\sm I)=A\mod \sm I\subseteq \cat{C}\sm I$.  Thus the domains of
$I$ are small in $\cat{C}$ relative to $(A\mod \sm _{A}(A\sm I))\cell $,
and so the domains of $A\sm I$ are small in $A\mod $ relative to $(A\mod
\sm _{A}(A\sm I))\cell $.  Similarly, the domains of $A\sm J$ are small
in $A\mod $ relative to $(A\mod \sm _{A}(A\sm J))\cell $.  And the maps
of $(A\mod \sm _{A}(A\sm J))\cell $ are in particular maps of
$(\cat{C}\sm J)\cell $, so are weak equivalences. Thus the category
$A\mod $ satisfies the same conditions as does $\cat{C}$ (except for the
second half of the definition of a monoidal model category), and so we
may as well assume that $A=S$.

It is well-known that the category $S\alg $ is bicomplete, and the two
out of three and retract axioms are immediate consequences of the
definitions.  The cofibration-trivial fibration half of the lifting
axiom is also immediate.  For the factorization axioms, we use the sets
$T(I)$ and $T(J)$.  Adjointness guarantees that the cofibrations in
$S\alg $ are the elements of $T(I)\cof $, and the trivial fibrations are
the elements of $T(I)\inj $.  Similarly, the fibrations are the elements
of $T(J)\inj $.  To understand the maps of $T(I)\cell $, we need to
understand the pushout in $S\alg $ of a map $T(g)\mathcolon
T(K)\xrightarrow{}T(L)$ through a map $T(K)\xrightarrow{}X$.  This
pushout is described in~\cite[Lemma 5.2]{schwede-shipley-monoids} as a
countable composition of maps $P_{i}\xrightarrow{h_{i}}P_{i+1}$.  Each
map $h_{i}$ is a pushout in $\cat{C}$ of $X^{\sm (i+1)}\sm g_{i}$, where
$g_{i}$ is very similar to the inclusion of the fat wedge into the
product.  That is, the target of $g_{i}$ is $L^{\sm (i)}$, and the
source is analogous to the subset of the smash product where at least
one term is in $K$.  In any case, one can see from the fact that
$\cat{C}$ is monoidal that $g_{i}$ is a cofibration if $g\in I$, and is
a trivial cofibration if $g\in J$.  It follows that the maps of
$T(I)\cell $ are in $(\cat{C}\sm I\cof )\cell $, which is contained in
$(\cat{C}\sm I)\cof $.  Since the domains of $I$ are small relative to
$(\cat{C}\sm I)\cell $, they are also small relative to $(\cat{C}\sm
I)\cof $ by~\cite[Theorem 12.4.21]{hirschhorn}.  Thus, by adjointness,
the domains of $T(I)$ are small in $S\alg $ relative to $T(I)\cell $.
The small object argument then gives the desired factorization into a
cofibration followed by a trivial fibration.

A very similar argument shows that the domains of $J$ are small relative
to $T(J)\cell $, and also that the maps of $T(J)\cell $, since they are
in $(\cat{C}\sm J)\cof $, are weak equivalences.  Thus the small object
argument applied to $T(J)$ gives the desired factorization into a
trivial cofibration followed by a fibration.  The proof that trivial
cofibrations have the \llp fibrations then uses the Retract Argument, as
in the proof of Theorem~\ref{thm-cofibrant-modules}.  
\end{proof}

\begin{example}\label{ex-top-monoids}
For an example where Theorem~\ref{thm-monoids-yes} applies but the
simpler version of~\cite{schwede-shipley-monoids} does not, let
$\cat{C}$ be the category of compactly generated spaces.  We have
already seen that the monoid axiom holds here (see
Lemma~\ref{lem-monoid-axiom-fibrant}).  One can also verify that the
elements of $(\cat{C}\sm I)\cell $ are closed inclusions, using the
results of~\cite[Appendix A]{lewis-thesis}.  It follows that the domains
of $I$ are small relative to $(\cat{C}\sm I)\cell $, and also that the
domains of $J$ are small relative to $(\cat{C}\sm J)\cell $.  Thus we do
get a model category of monoids of compactly generated spaces. On the
other hand, if we let $\cat{C}$ be the category of $k$-spaces, the
monoid axiom holds, but the necessary smallness conditions do not, so
far as we know.  So we do not get a model category of monoids of
$k$-spaces.
\end{example}

It is interesting to note that, in the situation of
Theorem~\ref{thm-monoids-yes}, the category $A\alg $ is a model category
even though the category $A\mod $ may not be a monoidal model category
unless either $A$ or $S$ is cofibrant.  

In any case, we are really interested in the case where the hypotheses
of Theorem~\ref{thm-monoids-yes} do not hold.  In this case, we obtain
the following result.  

\begin{theorem}\label{thm-monoids-almost}
Suppose $\cat{C}$ is a cofibrantly generated symmetric monoidal model
category.  Suppose that $A$ is either $S$ or a commutative monoid which
is cofibrant in $\cat{C}$.  Then the category $A\alg $ is almost a
model category, in the following precise sense.
\begin{enumerate}
\item $A\alg $ is bicomplete and the two out of three and retract
axioms hold. 
\item Cofibrations have the \llp trivial fibrations, and trivial
cofibrations whose source is cofibrant in $A\mod $ have the \llp
fibrations.
\item Every map whose source is cofibrant in $A\mod $ can be
functorially factored into a cofibration followed by a trivial
fibration, and also can be functorially factored into a trivial
cofibration followed by a fibration.
\end{enumerate}
Furthermore, cofibrations whose source is cofibrant in $A\mod $ are
cofibrations in $A\mod $, and fibrations and trivial fibrations are
closed under pullbacks.
\end{theorem}

\begin{proof}
The category $A\alg $ is the category of monoids in $A\mod $, which
itself is a cofibrantly generated monoidal model category, by
Proposition~\ref{prop-boxprod-R}.  Thus we can assume that $A=S$.  We
have already seen that bicompleteness and the retract and two out of
three axioms hold.  The lifting axiom for cofibrations and trivial
fibrations holds by definition.  As before, adjointness implies that the
trivial fibrations form the class $T(I)\inj $, so the cofibrations form
the class $T(I)\cof $.  The fibrations form the class $T(J)\inj $, so
the elements of $T(J)\cof $ have the \llp fibrations.  Recall from the
proof of Theorem~\ref{thm-monoids-yes} that the pushout in $S\alg $ of a
map $T(K)\xrightarrow{T(g)}T(L)$ through a map $T(K)\xrightarrow{}X$ of
monoids is a countable composition of maps
$P_{i}\xrightarrow{f_{i}}P_{i+1}$, where $f_{i}$ is the pushout in
$\cat{C}$ of a map $X^{\sm (i+1)}\sm g_{i}$.  The map $g_{i}$ is, as we
have said, a cofibration if $g$ is so, and a trivial cofibration if $g$
is so.  Thus, if $X$ is cofibrant in $\cat{C}$, a pushout in
$\monoids{C}$ of a map of $T(I)$ through a map to $X$ is a cofibration
in $\cat{C}$ and a pushout of a map of $T(J)$ through a map to $X$ is a
trivial cofibration in $\cat{C}$.  Since the domains of $I$ are small
relative to cofibrations in $\cat{C}$, it follows by adjointness that
the domains of $T(I)$ are small relative to transfinite compositions of
pushouts of maps of $T(I)$ as long as the initial stage is cofibrant in
$\cat{C}$.  The small object argument then allows us to functorially
factor any map in $S\alg $ whose source is cofibrant in $\cat{C}$
into a cofibration followed by a trivial fibration.  A similar argument,
with $T(J)$ replacing $T(I)$, allows us to functorially factor any map
in $S\alg $ whose source is cofibrant in $\cat{C}$ into a trivial
cofibration followed by a fibration.

Now, given any cofibration $f$ in $\monoids{C}$ whose source is
cofibrant, we can factor $f=pi$, where $i$ is a transfinite composition
of pushouts of maps of $T(I)$, and is therefore a cofibration in
$\cat{C}$, and $p$ is a trivial fibration.  The retract argument shows
that $f$ is a retract of $i$, and so is also a cofibration in $\cat{C}$.
Now suppose $f$ is a trivial cofibration, again with source cofibrant in
$\cat{C}$.  Then we can factor $f=qj$, where $q$ is a fibration of
monoids and $j$ is a transfinite composition of pushouts of $T(J)$.
Since $j$ is also a weak equivalence, so is $q$.  Thus $f$ has the \llp
$q$, and so $f$ is a retract of $j$.  Then $f$ has the \llp fibrations,
as desired.
\end{proof}

We should point out that it is not even clear that the cofibrations
$T(J)$ are weak equivalences in $\monoids{C}$.  This will be true if we
can choose the domains of $J$ to be cofibrant.

The information contained in Theorem~\ref{thm-monoids-almost} is enough
to carry out most of the usual constructions in model category theory,
at least as long as we have some element of $A\alg $ which is cofibrant
in $A\mod $.  For this, we must require that the unit $S$ is cofibrant
in $\cat{C}$.  Then the initial object $A$ of $A\alg $ is cofibrant in
$A\mod $.  Hence there is a cofibrant replacement functor $Q$ in
$A\alg $ and a natural trivial fibration $QX\xrightarrow{}X$.
There is a fibrant replacement functor $R$ defined on $A$-algebras which are
cofibrant in $A\mod $, in particular on cofibrant $A$-algebras, and a
natural trivial cofibration $X\xrightarrow{}RX$.  In particular, though
the composite functor $QR$ does not make sense in general, the composite
$RQ$ does, and we find that the quotient category of $A\alg $
obtained by inverting the weak equivalences is equivalent to the
quotient category of the cofibrant and fibrant objects.

One does have to be careful though.  We can not recognize trivial
cofibrations in $A\alg $ as the class of maps with the \llp 
fibrations.  But a map in $A\alg $ whose source is
cofibrant in $A\mod $ is a trivial cofibration if and only if it has
the \llp all fibrations.  Similarly, one does not know that a pushout of a
trivial cofibration $f$ in $A\alg $ through a map $g$ is a trivial
cofibration unless both the source of $f$ and the target of $g$ are
cofibrant in $A\mod $.  Cylinder and path objects need not exist unless
the object in question is cofibrant in $A\mod $.  

Nevertheless, we can follow the usual definitions in the standard
construction of the homotopy category of a model category.  See for
example~\cite[Chapter 1]{hovey-model}, \cite{dwyer-spalinski},
or~\cite[Chapters 8 and 9]{hirschhorn}.  In order for the notion of left
homotopy to have any content, one must assume that the sources of one's
maps are cofibrant in $A\mod $.  Similarly, for right homotopy, one
must assume that the target is cofibrant in $A\mod $.  It is not clear
that right homotopy is an equivalence relation if the target is
cofibrant in $A\mod $ and fibrant, as one would expect.  One must first
prove that the notions of left and right homotopy coincide if the source
is cofibrant and the target is cofibrant and fibrant.  Then, since left
homotopy is an equivalence relation in that situation, so is right
homotopy.  

In the end, we obtain the following theorem.

\begin{theorem}\label{thm-monoids-homotopy}
Suppose $\cat{C}$ is a cofibrantly generated symmetric monoidal model
category, where the unit $S$ is cofibrant.  Suppose $A$ is a commutative
monoid which is cofibrant in $\cat{C}$.  Let $\ho A\alg $ denote the
category obtained from $A\alg $ by formally inverting the weak
equivalences.  Then there is an equivalence of categories
$(A\alg)_{cf}/\sim \xrightarrow{}\ho A\alg $, where
$(A\alg )_{cf}$ is the full subcategory of cofibrant and fibrant
$A$-algebras, and $\sim $ denotes the homotopy equivalence relation.  In
particular, $\ho A\alg $ exists.  A map in $A\alg $ is an
isomorphism in $\ho A\alg $ if and only if it is a weak
equivalence.
\end{theorem}

We now show that the homotopy category of monoids is homotopy invariant. 
The following theorem is analogous to
Theorem~\ref{thm-maps-cofibrant-monoids}.  

\begin{theorem}\label{thm-maps-cofibrant-algebras}
Suppose $\cat{C}$ is a cofibrantly generated symmetric monoidal model
category such that the unit $S$ is cofibrant and the domains of the
generating cofibrations can be taken to be cofibrant.  Suppose
$f\mathcolon A\xrightarrow{}A'$ is a weak equivalence of commutative
monoids which are cofibrant in $\cat{C}$.  Then the derived functors of
induction and restriction define an adjoint equivalence of categories
$\ho A\alg \xrightarrow{}\ho A'\alg $.  
\end{theorem}

\begin{proof}
One can easily check that induction which takes $M$ to $A'\sm _{A}M$,
defines a functor from $A$-algebras to $A'$-algebras, which is left
adjoint to the restriction functor.  The restriction functor obviously
preserves fibrations and trivial fibrations, and reflects weak
equivalences between fibrant objects.  It follows that induction
preserves cofibrations, and those trivial cofibrations whose source is
cofibrant in $A\mod $.  One can then easily check that the functor which
takes $M$ to $A'\sm _{A}QM$ is the total left derived functor of
induction, and that its right adjoint is the functor which takes $M$ to
$RQM$.  The situation is very similar to the usual Quillen functor
formalism, it is just that we need to apply $Q$ before applying the
fibrant replacement functor $R$ since $R$ is not globally defined.

The same argument as in Lemma~\ref{lem-Quillen-equiv} applies, so we are
reduced to showing that the map $M\xrightarrow{}A'\sm _{A}M$ is a weak
equivalence for all cofibrant $A$-algebras $M$.  But since every
cofibrant $A$-algebra is also a cofibrant $A$-module,
Theorem~\ref{thm-maps-cofibrant-monoids} completes the proof.  
\end{proof}

Then the following theorem is analogous to
Theorem~\ref{thm-modules-category}.  

\begin{theorem}\label{thm-invariance-monoids}
Suppose $F\mathcolon \cat{C}\xrightarrow{}\cat{D}$ is a monoidal Quillen
equivalence of cofibrantly generated symmetric monoidal model
categories.  Suppose as well that $S$ is cofibrant in $\cat{C}$, and $A$
is a commutative monoid which is cofibrant in $\cat{C}$.  Then $F$
induces an equivalence of categories $LF\mathcolon \ho
A\alg \xrightarrow{}\ho FA\alg $.
\end{theorem}

\begin{proof}
Let $U$ denote the right adjoint of $F$.  It is clear that $F$ defines a
functor from $A$-algebras to $FA$-algebras.  Indeed, $F$ defines a
monoidal functor from $A$-modules to $FA$-modules.  To see that $U$ defines a
functor going the other way, note that we have a natural map $UX\sm_{A}
UY\xrightarrow{}U(X\sm_{FA} Y)$ adjoint to the composite
\[
F(UX\sm_{A} UY) \cong FUX \sm_{FA} FUY \xrightarrow{\varepsilon \sm
\varepsilon } X\sm_{FA} Y
\]
Thus, if $X$ is an $A$-algebra, so is $UX$; the unit of $UX$ is adjoint to the
unit of $X$ using the isomorphism.  One can easily check
that $F\mathcolon A\alg \xrightarrow{}FA\alg $ is left adjoint
to $U\mathcolon FA\alg \xrightarrow{}A\alg $.  

The functor $U$ preserves fibrations and trivial fibrations of monoids.
It follows that $F$ preserves cofibrations and those trivial
cofibrations with source cofibrant in $A\mod $.  We can then define
$(LF)X=F(QX)$ as usual.  We have to define $(RU)X=U(RQX)$ however, since
$RX$ does not make sense in general.  The usual argument shows that
$LF\mathcolon \ho A\alg \xrightarrow{}\ho FA\alg $ is left
adjoint to $RU$.  

The functor $U$ reflects weak equivalences between fibrant objects.  A
similar argument as in Lemma~\ref{lem-Quillen-equiv} implies that we
need only check that the map $X\xrightarrow{}URFX$ is a weak equivalence
for cofibrant $A$-algebras $X$.  Note that $FX$ is still cofibrant, so
we do not need to consider $URQFX$, though of course we could do so
without difficulty.  Here $R$ is a fibrant replacement functor in
$FA\alg $.  If we let $L$ denote a fibrant replacement functor in
$FA\mod $, then we know from Theorem~\ref{thm-modules-category} that
$X\xrightarrow{}ULFX$ is a weak equivalence, since $X$ is also cofibrant
as an $A$-module.  A simple lifting argument shows that there is a weak
equivalence $LFX\xrightarrow{}RFX$ in $FA\mod$.  Since $U$ preserves
weak equivalences between fibrant objects, this gives a weak equivalence
$ULFX\xrightarrow{}URFX$.  Hence $X\xrightarrow{}URFX$ is a weak
equivalence as well.
\end{proof}

\begin{example}\label{ex-spectra}
For example, we can take $\cat{C}$ to be the symmetric monoidal model
category of simplicial symmetric spectra~\cite{hovey-shipley-smith}.
Here the monoid axiom holds and everything is small relative to the
whole category, so we get a model category of $S$-algebras.  We take
$\cat{D}$ to be the symmetric monoidal model category of topological
symmetric spectra.  We do not know whether the monoid axiom holds here.
Nonetheless, the geometric realization is a monoidal Quillen equivalence
$\cat{C}\xrightarrow{}\cat{D}$, so defines an equivalence of categories
$\ho S\alg \xrightarrow{}\ho S\alg $ between the homotopy
categories of $S$-algebras.  
\end{example}

Note that the homotopy category of an arbitrary model category has a
closed action by the homotopy category of simplicial sets~\cite[Chapter
5]{hovey-model}.  This result will still hold for $\ho R\alg $; one
can use the same approach as in~\cite[Chapter 5]{hovey-model}, replacing
$R$ by $RQ$ everywhere.  In particular, there are mapping spaces of
monoids.  The equivalence of Theorem~\ref{thm-invariance-monoids} will
preserve these mapping spaces.


\begin{thebibliography}{EKMM97}

\bibitem[DHK]{kan-model}
W.~G. Dwyer, P.~S. Hirschhorn, and D.~M. Kan, \emph{Model categories and
  general abstract homotopy theory}, in preparation.

\bibitem[DS95]{dwyer-spalinski}
W.~G. Dwyer and J.~Spalinski, \emph{Homotopy theories and model categories},
  Handbook of algebraic topology (Amsterdam), North-Holland, Amsterdam, 1995,
  pp.~73--126.

\bibitem[EKMM97]{elmendorf-kriz-mandell-may}
A.~D. Elmendorf, I.~Kriz, M.~A. Mandell, and J.~P. May, \emph{Rings, modules,
  and algebras in stable homotopy theory. {W}ith an appendix by {M}. {C}ole},
  Mathematical Surveys and Monographs, vol.~47, American Mathematical Society,
  Providence, RI, 1997, xii+249 pp.

\bibitem[Hir97]{hirschhorn}
P.~S. Hirschhorn, \emph{Localization, cellularization, and homotopy colimits},
  preprint, 1997.

\bibitem[Hov97]{hovey-model}
Mark Hovey, \emph{Model categories}, preprint (x+194 pages), 1997.

\bibitem[Hov98]{hovey-stable-model}
Mark Hovey, \emph{Stabilization of model categories}, preprint, 1998.

\bibitem[HSS98]{hovey-shipley-smith}
Mark Hovey, Brooke~E. Shipley, and Jeffrey~H. Smith, \emph{Symmetric spectra},
  preprint, 1998.

\bibitem[Lew78]{lewis-thesis}
L.~G. Lewis, Jr., \emph{The stable category and generalized {T}hom spectra},
  Ph.D. thesis, University of Chicago, 1978.

\bibitem[Qui67]{quillen-htpy}
Daniel~G. Quillen, \emph{Homotopical algebra}, Lecture Notes in Mathematics,
  vol.~43, Springer-Verlag, 1967.

\bibitem[SS97]{schwede-shipley-monoids}
Stefan Schwede and Brooke Shipley, \emph{Algebras and modules in monoidal model
  categories}, preprint, 1997.

\end{thebibliography}

\providecommand{\bysame}{\leavevmode\hbox to3em{\hrulefill}\thinspace}

\end{document}